\documentclass[11pt]{amsart}
\usepackage{amsmath,amssymb,amsthm,upref,graphicx,mathrsfs, enumerate}
\usepackage{esint}
\usepackage{textcomp} 
\usepackage{array} 
\usepackage{textcomp} 
\usepackage[
  colorlinks=true,
  linkcolor=blue,
  citecolor=blue,
  urlcolor=blue,
  pagebackref]{hyperref}

\numberwithin{equation}{section}
\allowdisplaybreaks

\newtheorem{theorem}{Theorem}[section]
\newtheorem{prop}[theorem]{Proposition}

\newtheorem{cor}[theorem]{Corollary}

\theoremstyle{definition}
\newtheorem{definition}[theorem]{Definition}

\newtheorem{remark}[theorem]{Remark}
\theoremstyle{theorem}

\def\XXint#1#2#3{{\setbox0=\hbox{$#1{#2#3}{\int}$}
\vcenter{\hbox{$#2#3$}}\kern-.5\wd0}}

\begin{document}

\title[Relations between Schramm spaces and generalized Wiener classes]{Relations between Schramm spaces and generalized Wiener classes}

\author[Milad Moazami Goodarzi, Mahdi Hormozi and Nacima Memi\'{c}]{Milad Moazami Goodarzi, Mahdi Hormozi and Nacima Memi\'{c}}

\address[M. Moazami Goodarzi]{Department of Mathematics, Faculty of Sciences, Shiraz University, Shiraz 71454, Iran}
\email{milad.moazami@gmail.com}

\address[M. Hormozi]{Department of Mathematical Sciences, Division of Mathematics,
University of Gothenburg, Gothenburg 41296, Sweden \& Department of Mathematics, Faculty of Sciences, Shiraz University, Shiraz 71454, Iran}
\email{me.hormozi@gmail.com}

\address[N. Memi\'{c}]
{Department of Mathematics, Faculty of Natural Sciences and Mathematics, University of Sarajevo, Zmaja od Bosne 33-35, Sarajevo, Bosnia and Herzegovina}
\email{nacima.o@gmail.com}


\date{\today}


\begin{abstract}
We give necessary and sufficient conditions for the embeddings $\Lambda\text{BV}^{(p)}\subseteq \Gamma\text{BV}^{(q_n\uparrow q)}$ and $\Phi\text{BV}\subseteq\text{BV}^{(q_n\uparrow q)}$. As a
consequence, a number of results in the literature, including a fundamental theorem of Perlman and Waterman, are simultaneously extended.
\end{abstract}

\maketitle


\section{\bf{Introduction and main results}}\label{Sect:1}

Let $\Lambda=\{\lambda_j\}_{j=1}^\infty$ be a nondecreasing sequence of positive numbers such that $\sum_{j=1}^\infty\frac1{\lambda_j}=\infty\:$. Following
\cite{Appell}, we call $\Lambda$ a Waterman sequence.
Let $\Phi=\{\phi_j\}_{j=1}^\infty$ be a sequence of increasing convex functions on $[0,\infty)$ with $\phi_j(0)=0$.
We say that $\Phi$ is a Schramm sequence if $0<\phi_{j+1}(x)\leq\phi_j(x)$ for all $j$ and $\sum_{j=1}^\infty\phi_j(x)=\infty$ for all $x>0$. This terminology is used throughout.

We begin by recalling two generalizations of the concept of bounded variation which are central to our work.
\begin{definition}
A real-valued function $f$ on $[a,b]$ is said to be of $\Phi$-bounded
variation if
$$
V_\Phi(f)=V_\Phi(f;[a,b])=\sup\sum_{j=1}^n\phi_j(|f(I_j)|)<\infty,
$$
where the supremum is taken over all finite collections $\{I_j\}_{j=1}^n$ of nonoverlapping subintervals of $[a,b]$ and $f(I_j)=f(\sup I_j)-f(\inf I_j)$.
We denote by $\Phi\text{BV}$ the linear space of all functions $f$ such that $cf$ is of $\Phi$-bounded variation for some $c>0$.
\end{definition}
If for every $f\in\Phi\text{BV}$, we define
$$
\|f\|:=|f(a)|+\inf\{c>0:V_\Phi(f/c)\leq1\},
$$
then it is easily seen that $\|\cdot\|$ is a norm, and $\Phi\text{BV}$ endowed with this norm turns into a Banach space. The space $\Phi\text{BV}$ is introduced in Schramm's
paper \cite{SW2}. For more information about $\Phi\text{BV}$, the reader is referred to \cite{Appell}.

If $\phi$ is a strictly increasing convex function on $[0,\infty)$ with $\phi(0)=0$, and if $\Lambda=\{\lambda_j\}_{j=1}^\infty$ is
a Waterman sequence, by taking $\phi_j(x)=\phi(x)/\lambda_j$ for all $j$, we get the class $\phi\Lambda\text{BV}$ of functions of $\phi\Lambda$-bounded variation. This
class was introduced by Schramm and Waterman in \cite{SW1} (see also \cite{SW} and \cite{L}). More specifically, if $\phi(x)=x^p$ ($p\geq1$), we get the Waterman-Shiba class
$\Lambda\text{BV}^{(p)}$, which was introduced by Shiba in \cite{S}. When $p=1$, we obtain the well-known Waterman class $\Lambda\text{BV}$.

In the case $\lambda_j=1$ for all $j$, we obtain the class $\phi\text{BV}$ of functions of $\phi$-bounded variation introduced by Young \cite{You}. More specifically,
when $\phi(x)=x^p$ ($p\geq1$), we obtain the Wiener class $\text{BV}_p$ (see \cite{Wi}), and taking $p=1$, we have the well-known Jordan class $\text{BV}$.

\begin{remark}\label{simple}
One can easily observe that functions of $\Phi$-bounded variation are bounded and can only have simple discontinuities (countably many of them, indeed). The
class $\Phi\text{BV}$ has many applications in Fourier analysis as well as in treating topics such as convergence, summability, etc. (see \cite{Wi,You,Wt1,Wt2,Wt3,Perlman,SW2}).
\end{remark}

\begin{definition}
Let $\{q_n\}_{n=1}^\infty$ and $\{\delta_n\}_{n=1}^\infty$ be sequences of positive real numbers such that $1\leq q_n\uparrow q\leq\infty$ and $2\leq\delta_n\uparrow\infty$.
A real-valued function $f$ on $[a,b]$ is said to be of $q_n$-$\Lambda$-bounded variation if
$$
V_\Lambda(f)=V_\Lambda(f;q_n\uparrow q;\delta):=\sup_{n\geq 1}\sup_{\{I_j\}}\Big(\sum_{j=1}^{s} \frac{|f(I_j)|^{q_n}}{\lambda_j}\Big)^{\tfrac{1}{q_n}}<\infty,
$$
where the $\{I_j\}_{j=1}^s$ are collections of nonoverlapping subintervals of $[a,b]$ such that  $\inf_{j} |I_j|\geq \frac{b-a}{\delta_n}$. The class of functions
of $q_n$-$\Lambda$-bounded variation is denoted by $\Lambda\text{BV}^{(q_n\uparrow q)}$ ($=\Lambda\text{BV}_\delta^{(q_n\uparrow q)}$). In the sequel, we suppose that
$[a,b]=[0,1]$.
\end{definition}
The class $\Lambda\text{BV}^{(q_n\uparrow q)}$ was introduced by Vyas in \cite{Vyas}. When $\lambda_j=1$ for all $j$ and $\delta_n=2^n$ for all $n$, we get the class
$\text{BV}^{(q_n\uparrow q)}$---inroduced by Kita and Yoneda (see \cite{KY})---which in turn recedes to the Wiener class $\text{BV}_q$, when $q_n=q$ for all $n$.

\bigskip
A natural and important problem is to determine relations between the above-mentioned classes; see \cite{Wt1}, \cite{Perlman}, \cite{Ciem}, \cite{KY},
\cite{G}, \cite{Pierce}, \cite{HPR} and \cite{Wang2} for some results in this direction. In particular, Perlman and Waterman found the fundamental characterization of embeddings between
$\Lambda\text{BV}$ classes in \cite{Perlman}. Ge and Wang characterized the embeddings $\Lambda\text{BV}\subseteq\phi\text{BV}$ and $\phi\text{BV}\subseteq\Lambda\text{BV}$
(see \cite{Wang2}). It was shown by Kita and Yoneda in \cite{KY} that the embedding
$\text{BV}_p\subseteq\text{BV}^{(p_n\uparrow\infty)}$ is both automatic and strict for all $1\leq p<\infty.$
Furthermore, Goginava characterized the embedding $\Lambda\text{BV}\subseteq\text{BV}^{(q_n\uparrow\infty)}$ in \cite{G}, and a characterization of
the embedding $\Lambda\text{BV}^{(p)}\subseteq\text{BV}^{(q_n\uparrow q)}$ ($1\leq q\leq\infty$) was given by Hormozi, Prus-Wi\'{s}niowski and Rosengren in \cite{HPR}.
In this paper, we investigate the embeddings $\Lambda\text{BV}^{(p)}\subseteq \Gamma\text{BV}^{(q_n\uparrow q)}$ and $\Phi\text{BV}\subseteq\text{BV}^{(q_n\uparrow q)}$ ($1\leq q\leq\infty$). The problem as to when the reverse embeddings hold is also considered, which turns out to have a simple answer (see Remark \eqref{reverse}(ii)
below).

\bigskip
Throughout this paper, the letters $\Lambda$ and $\Gamma$ are reserved for a typical Waterman sequence. We associate to $\Lambda$ a function which we still denote by $\Lambda$
and define it as $\Lambda(r):=\sum_{j=1}^{[r]} \frac1{\lambda_j}$ for $r\geq 1$. The function $\Lambda(r)$ is clearly nondecreasing and
$\Lambda(r)\rightarrow\infty$ as $r\rightarrow\infty$. Our first main result reads as follows.
\begin{theorem}\label{t1}
Let $1\leq p\leq q_n\uparrow q\leq\infty$. Then, a necessary and sufficient condition for the embedding $\Lambda\emph{BV}^{(p)}\subset \Gamma\emph{BV}^{(q_n\uparrow q)}$ is
\vspace{.25cm}
\begin{equation}\label{th1}
\limsup_{n\rightarrow\infty}\max_{1\leq k\leq \delta_n}\Gamma(k)^{\frac{1}{q_n}}\Lambda(k)^{-\frac{1}{p}}< \infty.
\end{equation}
Moreover, if the hypothesis is replaced by the condition that $\{\Gamma(n)/\Lambda(n)\}_{n=1}^\infty$ be nondecreasing, then the conclusion of the theorem
still holds true.
\end{theorem}
\bigskip
An important consequence of Theorem \eqref{t1} is the following corollary, which is indeed a nontrivial extension of \cite[Theorem 3]{Perlman}.
\begin{cor}
Let $1\leq p\leq q<\infty$. Then, a necessary and sufficient condition for the embedding $\Lambda\emph{BV}^{(p)}\subseteq \Gamma\emph{BV}^{(q)}$ is
$$
\sup_{1\leq n<\infty}\frac{\Gamma(n)^{\frac{1}{q}}}{\Lambda(n)^{\frac{1}{p}}}< \infty.
$$
\end{cor}


\begin{cor}\emph{(\cite[Theorem 1]{HPR})}
Let $1\leq p<\infty$. Then, a necessary and sufficient condition for the embedding $\Lambda\emph{BV}^{(p)}\subseteq\emph{BV}^{(q_n\uparrow q)}$ is
$$
\limsup_{n\rightarrow\infty}\max_{1\leq k\leq \delta_n}\frac{k^{\frac{1}{q_n}}}{\left(\sum_{i=1}^k \frac1{\lambda_i}\right)^{\frac{1}{p}}}< \infty.
$$
\end{cor}
\bigskip
Next corollary extends \cite[Lemma 2.1]{KY}.

\begin{cor}
Let $1<q\leq\infty$. Then, We have
$$
\bigcup_{1\leq p<q}\Lambda\emph{BV}^{(p)}\subseteq\Lambda\emph{BV}^{(q_n\uparrow q)}.
$$
\end{cor}
\bigskip
If $\Phi=\{\phi_j\}_{j=1}^\infty$ is a Schramm sequence, we define $\Phi_k(x):=\sum_{j=1}^k\phi_j(x)$ for $x\geq0$. Then $\Phi_k(x)$ is clearly an increasing
convex function on $[0,\infty)$ such that $\Phi_k(0)=0$ and $\Phi_k(x)>0$ for $x>0$. Without loss of generality we assume that $\Phi_k(x)$ is strictly
increasing on $[0,\infty)$. Let $\Phi_k^{-1}(x)$ be the inverse function of $\Phi_k(x)$. Our next main result can be formulated as follows.

\begin{theorem}
\label{t2}
A necessary and sufficient condition for the embedding $\Phi\emph{BV}\subset\emph{BV}^{(q_n\uparrow q)}$ is
\vspace{.25cm}
\begin{equation}\label{th3}
\limsup_{n\rightarrow\infty}\max_{1\leq k\leq \delta_n}k^{\frac{1}{q_n}}\Phi_k^{-1}(1)< \infty.
\end{equation}
\end{theorem}

\begin{cor}
A necessary and sufficient condition for the embedding $\phi\Lambda\emph{BV}\subset\emph{BV}^{(q_n\uparrow q)}$ is
$$
\limsup_{n\rightarrow\infty}\max_{1\leq k\leq \delta_n}k^{\frac{1}{q_n}}\phi^{-1}(\Lambda(k)^{-1})< \infty.
$$
\end{cor}

\begin{remark}\label{reverse}
(i) When $\phi(x)=x^p$, $1\leq p<\infty$, Corollary 1.9 yields Corollary 1.6 as a special case.

(ii) By \cite[Theorem 3.3]{KY}, the class $\text{BV}^{(q_n\uparrow\infty)}$ always contains a function with nonsimple discontinuities. Since clearly
$\text{BV}^{(q_n\uparrow\infty)}\subseteq\Lambda\text{BV}^{(q_n\uparrow\infty)}$, this is also the case for the class $\Lambda\text{BV}^{(q_n\uparrow\infty)}$.
On the other hand, as pointed out in Remark \eqref{simple}, the functions in the classes $\Phi\text{BV}$ and $\Lambda\text{BV}^{(p)}$ can only have simple
discontinuities. Hence, the corresponding reverse embeddings can never happen.

\end{remark}


\section{\bf{An auxiliary inequality}}\label{Sect:2}

In this section we establish an inequality (see \eqref{ineq} below) which plays a crucial role in the sufficiency part of the
proof of Theorem \eqref{t1}. Also some applications of it are presented in Corollary \eqref{ap1} and Remark \eqref{ap3}. The following
proposition is indeed a generalization of \cite[Lemma]{kup}.

\begin{prop}
Let $1\leq q< \infty$ and $n\in\mathbb{N}$. Then
\begin{equation}\label{ineq}
\Big(\sum_{j=1}^n x_j^qz_j\Big)^{\frac1{q}}\leq\sum_{j=1}^n x_jy_j\max_{1\leq k\leq n}\Big(\sum_{j=1}^k z_j\Big)^{\frac1{q}}\Big(\sum_{j=1}^k y_j\Big)^{-1},
\end{equation}
where $\{x_j\}$, $\{y_j\}$ and $\{z_j\}$ are positive nonincreasing sequences.
\end{prop}
\noindent\textbf{Proof.} Without loss of generality we may assume that $\sum_{j=1}^n x_jy_j=1$. With this in mind, it is enough to prove that the maximum value of
$\sum_{j=1}^n x_j^qz_j$ under above assumptions is
$$
\max_{1\leq k\leq n}\Big(\sum_{j=1}^k z_j\Big)\Big(\sum_{j=1}^k y_j\Big)^{-q}.
$$

We claim that the solution to this problem satisfies condition

\begin{equation}
\label{order}
x_1=x_2=\cdots=x_k>x_{k+1}=x_{k+2}=\cdots=x_n=0
\end{equation}
for some $1\leq k\leq n$. To prove our claim, we suppose
to the contrary that there exists a solution which does not satisfy condition \eqref{order}. Then for some $1\leq k\leq n$, we have $x_{k+1}>0$ and
$$
x_1=x_2=\cdots=x_k>x_{k+1}\geq x_{k+2}\geq\cdots\geq x_n\geq 0.
$$
Put
$$
A:=\sum_{j=1}^k x_jy_j, \ \ B:=\sum_{j=k+1}^n x_jy_j, \ \ C:=\frac{x_{k+1}}{x_k},
$$
and define
$$
A\eta(t)+Bt=1.
$$
Then the $n$-tuple
$$
(\eta(t)x_1,\eta(t)x_2,\cdots,\eta(t)x_k,tx_{k+1},\cdots,tx_n)
$$
satisfies conditions of the problem, whenever $0\leq t<1/AC+B$. Now define
$$
f(t):=\eta(t)^q\sum_{j=1}^k x_j^qz_j+t^q\sum_{j=k+1}^n x_j^qz_j
$$
and consider two possibilities:

1) If $q>1$ then
$$
f''(t)=q(q-1)\Big(\eta(t)^{q-2}(\eta'(t))^2\sum_{j=1}^k x_j^qz_j+t^{q-2}\sum_{j=k+1}^n x_j^qz_j\Big)
$$
and hence $f''(1)>0$ which in turn implies that $f$ has a local minimum at $t=1$. This is a contradiction.

2) If $q=1$ then $f(t)$ is linear. Consequently,
$$
A\sum_{j=k+1}^n x_j^qz_j-B\sum_{j=1}^k x_j^qz_j=0
$$
which implies that the problem has
a solution satisfying condition \eqref{order}. This completes the proof.\ \ \ \ $\square$

\bigskip
Let $f$ be a bounded function on $[0,1]$. The modulus of variation of $f$ is the sequence $\nu_f$ and is defined by
$$
\nu_f(n):=\sup\sum_{j=1}^n |f(I_j)|,
$$
where the supremum is taken over all finite collections $\{I_j\}_{j=1}^n$ of nonoverlapping subintervals of $[0,1]$. The modulus of variation of $f$ is
nondecreasing and concave. A sequence $\nu$ with such properties is called a modulus of variation. The symbol $V[v]$ denotes the class of all functions $f$
for which there exists a constant $C>0$ (depending on $f$) such that $\nu_f(n)/\nu(n)\leq C$ for all $n$ (see \cite{Cha1}).
The following corollary is an immediate consequence of inequality \eqref{ineq}.

\begin{cor}\cite[Theorem 1]{Avdis}\label{ap1} The following inclusion holds.
$$
\Lambda\emph{BV}\subseteq V[n\Lambda(n)^{-1}].
$$
\end{cor}
\noindent\textbf{Proof.}
Let $\{I_j\}_{j=1}^n$ be a collection of nonoverlapping subintervals of $[0,1]$. If $f\in\Lambda\text{BV}$, $q=1$, $x_j=|f(I_j)|$, $y_j=1/\lambda_j$ and
$z_j=1$, from \eqref{ineq} we obtain
$$
\sum_{j=1}^n |f(I_j)|\leq\sum_{j=1}^n\frac{|f(I_j)|}{\lambda_j}\max_{1\leq k\leq n} k\Lambda(k)^{-1}\leq V_\Lambda(f)n\Lambda(n)^{-1},
$$
which means that $f\in V[n\Lambda(n)^{-1}]$.\ \ \ \ $\square$

\begin{remark}\label{ap3}
Let $\Lambda=\{\lambda_j\}$ and $\Gamma=\{\gamma_j\}$ be Waterman sequeneces. As stated on page 181 of \cite{Frank}, Perlman and Waterman have shown, in the course of the proof of
\cite[Theorem 3]{Perlman}, that if there is a constant $C$ such that
$$
\sum_{j=1}^n \frac1{\gamma_j}\leq C\sum_{j=1}^n \frac1{\lambda_j} \ \ \ \text{for all}\ n,
$$
then, given any nonincreasing sequence $\{a_j\}$ of nonnegative numbers,
$$
\sum_{j=1}^n \frac{a_j}{\gamma_j}\leq C\sum_{j=1}^n \frac{a_j}{\lambda_j}.
$$
It is worth mentioning that one can easily see that this is a simple consequence of inequality \eqref{ineq} above.
\end{remark}


\section{\textbf{Proofs of main results}}

\noindent\textbf{Proof of Theorem \eqref{t1}.}
\textbf{Necessity}. We proceed by contraposition. If \eqref{th1} does not hold, using the fact that $\Gamma(r)\rightarrow\infty$ as $r \rightarrow\infty$, we
may, without loss of generality, assume that $\gamma_1=1$ and for each $n$
 \begin{equation}
\label{first}
\Gamma(\delta_n) \geq 2^{n+2},
\end{equation}
and
 \begin{equation}
\label{first2}
\Gamma(r_n)^{\frac{1}{q_n}}\Lambda(r_n)^{-\frac{1}{p}}>\ 2^{4n}
\end{equation}
for some integer $r_n$, $1\leq r_n\leq \delta_n.$

\noindent We are going to construct a function $f$ in $\Lambda\text{BV}^{(p)}$ that does not belong to $\Gamma\text{BV}^{(q_n\uparrow q)}.$ To this end,
let $s_n$ be the greatest integer such that $2s_n-1\leq 2^{-n}\Gamma(\delta_n)$ and put $t_n=\min\{r_n,s_n\}.$
We define a sequence of functions $\{f_n\}_{n=1}^\infty$ on $[0,1]$ as follows:
$$ f_n(x):=\begin{cases}
 2^{-n}\Lambda(r_n)^{-\frac{1}{p}} ~~~~~, ~~~~x\in [2^{-n}+\frac{2j-2}{\delta_n},2^{-n}+\frac{2j-1}{\delta_n}) ;~~~~~  1\leq j\leq t_n, \\[.1in]
0~~~~~\qquad\qquad \textmd{otherwise}.\\
\end{cases}
$$
The functions $f_n$, defined in this fashion, have disjoint supports and therefore $f(x):=\sum_{n=1}^{\infty}f_n(x)$ is a well-defined function on $[0,1]$.
In addition, we have
\begin{align*}
V_\Lambda(f)\leq\sum_{n=1}^{\infty}V_\Lambda(f_n)=\sum_{n=1}^{\infty} \Big(\sum_{j=1}^{2t_n} \frac{(2^{-n}\Lambda(r_n)^{-\frac{1}{p}})^p}{\lambda_j}\Big)^{\frac1{p}}\\[.1in]
&\hspace{-7.3cm}\leq\sum_{n=1}^{\infty} 2^{-n+1}\Big(\sum_{j=1}^{r_n}\frac{\Lambda(r_n)^{-1}}{\lambda_j}\Big)^{\frac1{p}}
=\sum_{n=1}^{\infty} 2^{-n+1}\Big(\frac{\Lambda(r_n)}{\Lambda(r_n)}\Big)^{\frac1{p}}<\infty,
\end{align*}
since the sequence $\{\Lambda(r_n)^{-1}\}_{n=1}^\infty$ is nonincreasing and $t_n\leq r_n$. This means that $f\in\Lambda\text{BV}^{(p)}.$

On the other hand, $f\notin\Gamma\text{BV}^{(q_n\uparrow q)}.$ To see this, note that the definition of $s_n$ implies $2(s_n+1)-1>2^{-n}\Gamma(\delta_n).$ Combining
this with \eqref{first}, we obtain $\Gamma(2s_n-1)\geq 2^{-n-1}\Gamma(\delta_n).$ Consequently, if $t_n=s_n$, then the preceding inequality means that
$$
\Gamma(2t_n-1)\geq 2^{-n-1}\Gamma(\delta_n)\ge 2^{-n-1}\Gamma(r_n),
$$
since $r_n\leq\delta_n.$ Also, if $t_n=r_n$, clearly $\:2t_n-1\ge r_n$ and hence
$\Gamma(2t_n-1)\geq \Gamma(r_n)$, since $\Gamma(r)$ is increasing. Thus, we have shown
\begin{equation}
\label{second}
\Gamma(2t_n-1)\ge 2^{-n-1}\Gamma(r_n), \hspace{.3in}\text{for all $n$.}
\end{equation}

Finally, the intervals
$$
I_j:=\left[2^{-n}+\tfrac{j-1}{\delta_n},2^{-n}+\tfrac{j}{\delta_n}\right], \hspace{.3in}j=1,...,2t_n-1,
$$
have  length $\frac{1}{\delta_n}$ for each $n$, and thus
\begin{align*}
V_\Gamma(f)\ge\Big(\sum_{j=1}^{2t_n-1} \frac{|f(I_j)|^{q_n}}{\gamma_j}\Big)^{\frac1{q_n}}
=\left(\Gamma(2t_n-1)(2^{-n}\Lambda(r_n)^{-\frac{1}{p}})^{q_n}\right)^{\frac1{q_n}}\\[.1in]
&\hspace{-5.7cm}\ge 2^{-n}\left(2^{-n-1}\Gamma(r_n)(\Lambda(r_n)^{-\frac{1}{p}})^{q_n}\right)^{\frac1{q_n}}\ge 2^n,
\end{align*}
where the last two inequalities are due to \eqref{second} and \eqref{first2}, respectively. As a result, $V_\Gamma(f)$ is not finite.

\textbf{Sufficiency}. Assume \eqref{th1} and let $f\in\Lambda\text{BV}^{(p)}$. Let $\{I_j\}_{j=1}^s$ be a nonoverlapping collection of subintervals of $[0,1]$ with $\inf|I_j|\geq 1/\delta_n$,
and let $q=q_n/p\geq1$, $x_j=|f(I_j)|^p$, $y_j=1/\lambda_j$, $z_j=1/\gamma_j$. By \cite[Theorem 368]{Hardy}, we may also assume that the $x_j$'s are arranged in
descending order. Now, we can apply \eqref{ineq} and get
\begin{align*}\Big(\sum_{j=1}^s \frac{|f(I_j)|^{q_n}}{\gamma_j}\Big)^{\frac{1}{q_n}}
&\leq\Big(\sum_{j=1}^s \frac{|f(I_j)|^p}{\lambda_j}\Big)^{\frac{1}{p}}\max_{1\leq k\leq s} \Gamma(k)^{\frac{1}{q_n}}\Lambda(k)^{-\frac{1}{p}}\\
&\leq\Big(\sum_{j=1}^s \frac{|f(I_j)|^p}{\lambda_j}\Big)^{\frac{1}{p}}\max_{1\leq k\leq \delta_n} \Gamma(k)^{\frac{1}{q_n}}\Lambda(k)^{-\frac{1}{p}},
\end{align*}
where the second inequality is a consequence of $s\leq \delta_n.$ Taking suprema over all collections $\{I_j\}_{j=1}^s$ as above, and over all $n$ yields
$$
V_\Gamma(f) \leq V_\Lambda(f) \sup_{n}\max_{1\leq k\leq \delta_n}\Gamma(k)^{\frac{1}{q_n}}\Lambda(k)^{-\frac{1}{p}}< \infty.
$$
Hence $f\in\Gamma\text{BV}^{(q_n\uparrow q)}$ and the first part of the theorem is proved.

To prove the second part, let us assume that $\{\Gamma(n)/\Lambda(n)\}_{n=1}^\infty$ is nondecreasing. Observe that the proof of necessity is
identical to that given in the first part. For sufficiency, note that the only case which needs to be justified is when $q_n<p$ for some $n$. If this is the
case, we first apply \eqref{ineq} with $q=1$ to obtain
\begin{equation}\label{q=1}
\sum_{j=1}^s \frac{|f(I_j)|^p}{\gamma_j}\leq\sum_{j=1}^s \frac{|f(I_j)|^p}{\lambda_j}\max_{1\leq k\leq s}\Gamma(k)\Lambda(k)^{-1}.
\end{equation}
Then an application of H\"older's inequality yields
\begin{align*}
\sum_{j=1}^s \frac{|f(I_j)|^{q_n}}{\gamma_j}=\sum_{j=1}^s\Big(\frac{|f(I_j)|^p}{\gamma_j}\Big)^{\frac{q_n}{p}}\gamma_j^{\frac{q_n}{p}-1}\\[.1in]
&\hspace{-4.5cm}\leq\Big(\sum_{j=1}^s \frac{|f(I_j)|^p}{\gamma_j}\Big)^{\frac{q_n}{p}}\Gamma(s)^{1-\frac{q_n}{p}}\\[.1in]
&\hspace{-4.5cm}\leq\Big(\sum_{j=1}^s \frac{|f(I_j)|^p}{\lambda_j}\Big)^{\frac{q_n}{p}}\Gamma(s)^{1-\frac{q_n}{p}}\max_{1\leq k\leq s} \Gamma(k)^{\frac{q_n}{p}}\Lambda(k)^{-\frac{q_n}{p}}\\[.1in]
&\hspace{-4.5cm}\leq\Big(\sum_{j=1}^s \frac{|f(I_j)|^p}{\lambda_j}\Big)^{\frac{q_n}{p}}\max_{1\leq k\leq \delta_n} \Gamma(k)\Lambda(k)^{-\frac{q_n}{p}},
\end{align*}
where the last two inequalities are due, respectively, to \eqref{q=1} and the fact that $\{\Gamma(n)/\Lambda(n)\}_{n=1}^\infty$
is nondecreasing.\ \ \ \ $\square$

\bigskip
\noindent\textbf{Proof of Theorem \eqref{t2}.}
\textbf{Necessity}. Suppose \eqref{th3} does not hold. Then, without loss of generality, we may assume that for each $n$
$$
\delta_n\geq 2^{n+2},
$$
and
 \begin{equation}
\label{e7}
r_n^{\frac{1}{q_n}}\Phi_{r_n}^{-1}(1)>\ 2^{4n}
\end{equation}
for some integer $r_n$, $1\leq r_n\leq \delta_n.$

\noindent We will now construct a function $f\in\Phi\text{BV}$ such that $f\notin\text{BV}^{(q_n\uparrow q)}$. To do so, let $s_n$ be the greatest integer such that
$2s_n-1\leq 2^{-n}\delta_n$, let $t_n=\min\{r_n,s_n\}$ and consider the sequence $\{f_n\}_{n=1}^\infty$  of functions on $[0,1]$ defined in the following way:
$$ f_n(x):=\begin{cases}
 2^{-n}\Phi_{r_n}^{-1}(1) ~~~~~, ~~~~x\in [2^{-n}+\frac{2j-2}{\delta_n},2^{-n}+\frac{2j-1}{\delta_n}) ;~~~~~  1\leq j\leq t_n, \\[.1in]
0~~~~~\qquad\qquad \textmd{otherwise}.\\
\end{cases}
$$
Since the $f_n$'s have disjoint supports, $f(x):=\sum_{n=1}^{\infty}f_n(x)$ is a well-defined function on $[0,1]$. Thus, using convexity of the $\Phi_{r_n}$'s we have
\begin{align*}
V_\Phi(f) & \le\ \sum_{n=1}^{\infty}V_\Phi(f_n)\
= \sum_{n=1}^{\infty} \sum_{j=1}^{2t_n} \phi_j(2^{-n}\Phi_{r_n}^{-1}(1))
=\sum_{n=1}^{\infty} \Phi_{2t_n}(2^{-n}\Phi_{r_n}^{-1}(1)) \\[.1in]
&\hspace{2.6cm}\leq \sum_{n=1}^{\infty} \Phi_{2r_n}(2^{-n}\Phi_{r_n}^{-1}(1))
\leq\sum_{n=1}^{\infty} 2\Phi_{r_n}(2^{-n}\Phi_{r_n}^{-1}(1))<\infty,
\end{align*}
that is,  $f \in \Phi\text{BV}$.

In conclusion, let us show that $f\notin\text{BV}^{(q_n\uparrow q)}.$ To this end, proceeding in the same way as in the proof of Theorem \eqref{t1}, we obtain
\begin{equation}
\label{third}
2t_n-1\ge 2^{-n-1}r_n, \hspace{.3in}\text{for all $n$.}
\end{equation}
\noindent Since for every n, all intervals
$$
I_j:=\left[2^{-n}+\tfrac{j-1}{\delta_n},2^{-n}+\tfrac{j}{\delta_n}\right], \hspace{.3in}j=1,...,2t_n-1,
$$
have  length $\frac{1}{\delta_n}$, we get
\begin{align*}
V(f;q_n\uparrow q,\delta)\ge\Big(\sum_{j=1}^{2t_n-1}|f(I_j)|^{q_n}\Big)^{\frac1{q_n}}
=\Big((2t_n-1)(2^{-n}\Phi_{r_n}^{-1}(1))^{q_n}\Big)^{\frac1{q_n}}\\[.1in]
&\hspace{-5.2cm}\ge 2^{-n}\Big(2^{-n-1}r_n(\Phi_{r_n}^{-1}(1))^{q_n}\Big)^{\frac1{q_n}}
\ge 2^n,
\end{align*}
where the last two inequalities are results of \eqref{third} and \eqref{e7}, respectively. Therefore, $f\notin\text{BV}^{(q_n\uparrow q)}$.

\textbf{Sufficiency}. Let $f\in\Phi\text{BV}$. To show that $f\in\text{BV}^{(q_n\uparrow q)}$, it suffices to prove the inequality
\begin{equation}\label{final}
V(f;q_n\uparrow q;\delta)\leq C\sup_{n}\max_{1\leq k\leq \delta_n} k^{\frac{1}{q_n}}\Phi_k^{-1}(1),
\end{equation}
where $C$ is a positive constant depending solely on $f$.

\noindent In the course of the proof of Theorem 2.1 in \cite{Wu}, the author proceeds to estimate $(\sum_{j=1}^n x_j^q)^{\frac{1}{q}}$ under the restriction
$$
\sum_{j=1}^n \phi_j(x_{\tau(j)})\leq V_\Phi(f),
$$
where the $x_j$'s are arranged in descending order and $\tau$ is any permutation of $n$ letters. Using Wang's approach in \cite{Wang3}, he finds the following:
\begin{equation}
\label{Wu's estimate}
\Big(\sum_{j=1}^n x_j^q\Big)^{\frac{1}{q}}\leq 16\max_{1\leq k\leq n} k^{\frac1{q}}\Phi_k^{-1}(V_\Phi(f)).
\end{equation}

\noindent To prove \eqref{final}, consider a non-overlapping collection $\{I_j\}_{j=1}^s$ of subintervals of $[0,1]$
with $\inf|I_j|\geq 1/\delta_n$. If we put $q=q_n$,
$x_j=|f(I_j)|$, and if the $x_j$'s are rearranged in descending order, then we may apply \eqref{Wu's estimate} to obtain
\begin{align*}\Big(\sum_{j=1}^s |f(I_j)|^{q_n}\Big)^{\frac{1}{q_n}}
&\leq 16\max_{1\leq k\leq s} k^{\frac{1}{q_n}}\Phi_k^{-1}(V_\Phi(f))\\[.1in]
&\leq 16\max_{1\leq k\leq \delta_n} k^{\frac{1}{q_n}}\Phi_k^{-1}(V_\Phi(f)).\end{align*}
Taking suprema and using concavity of the $\Phi_k^{-1}$'s yields \eqref{final} with $C=16(1+V_\Phi(f))$.\ \ \ \ $\square$

\bigskip
\textbf{Acknowledgement.} The authors would like to thank Professor G.H. Esslamzadeh for kindly reading the manuscript of this paper and making valuable
remarks. The second author is supported by a grant from Iran's National Elites Foundation.


\begin{thebibliography}{EEE}
\bibitem{Appell}
J. Appell, J. Bana\'{s}, N. Merentes, Bounded Variation and Around, De Gruyter Ser. Nonlinear Anal. Appl., vol. 17, Walter de Gruyter, Berlin, 2013.

\bibitem{Avdis}
M. Avdispahi\'{c}, On the classes $\Lambda\text{BV}$ and $V[\nu]$, Proc. Amer. Math. Soc. 95 (2) (1985) 230--234.



\bibitem{Cha1}
Z.A. Chanturiya, The modulus of variation of a function and its application in the theory of Fourier series, Soviet. Math. Dokl. 15 (1974) 67--71.

\bibitem{Ciem}
J. Ciemnoczo\l owski, W. Orlicz, Inclusion theorems for classes of functions of generalized bounded variation, Comment. Math. 24 (1984) 181--194.


\bibitem{Wang2}
Y. Ge, H. Wang, Relationships between $\Phi\text{BV}$ and $\Lambda\text{BV}$, J. Math. Anal. Appl. 430 (2015) 1065--1073.

\bibitem{G}
U. Goginava, Relations between $\Lambda\text{BV}$ and $\text{BV}(p(n)\uparrow\infty)$ classes of functions, Acta Math. Hungar. 101 (4) (2003) 263-–272.

\bibitem{Hardy}
G. Hardy, J.E. Littlewood, G. Polya, Inequalities, 2nd edn., Cambridge University Press, Cambridge, 1952.

\bibitem{HPR}
M. Hormozi, and F. Prus-Wi\'{s}niowski, H. Rosengren, Inclusions of Waterman--Shiba spaces into generalized Wiener classes, J. Math. Anal. Appl. 419 (2014) 428--432.

\bibitem{KY}
H. Kita, K. Yoneda, A generalization of bounded variation, Acta Math. Hungar. 56 (3-4) (1990) 229--238.

\bibitem{kup}
Y.E. Kuprikov, Moduli of continuity of functions from Waterman classes, Moscow Univ. Math. Bull. 52 (5) (1997) 46--49.


\bibitem{L}
L. Leindler, A note on embedding of classes $H^\omega$, Anal. Math. 27 (2001) 71--76.


\bibitem{Perlman}
S. Perlman, D. Waterman, Some remarks on functions of $\Lambda$-bounded variation, Proc. Amer. Math. Soc. 74 (1979) 113--118.

\bibitem{Pierce}
P.B. Pierce, D.J. Velleman, Some generalization of the notion of bounded variation, Amer. Math. Monthly 113 (10) (2006) 897--904.

\bibitem{Frank}
F. Prus-Wi\'{s}niowski, Functions of bounded $\Lambda$-variation, In: Topics in Classical Analysis and Applications in Honor of Daniel Waterman, World
Sci. Publ. 173--190, Hackensack, 2008.


\bibitem{SW2}
M. Schramm,  Functions of $\Phi$-bounded variation and Riemann--Stieltjes integration, Trans. Amer. Math. Soc. 267 (1) (1985) 49--63.

\bibitem{SW1}
M. Schramm, D. Waterman, On the magnitude of Fourier coefficients, Proc. Amer. Math. Soc. 85 (1982) 407--410.

\bibitem{SW}
M. Schramm, D. Waterman, Absolute convergence of Fourier series of functions of $\Lambda\text{BV}^{(p)}$ and $\varphi\Lambda\text{BV}$, Acta Math. Hungar. 40
(3--4) (1982) 273--276.

\bibitem{S}
M. Shiba, On the absolute convergence of Fourier series of functions of class $\Lambda\text{BV}^{(p)}$, Sci. Rep. Fac. Ed. Fukushima Univ. 30 (1980) 7--10.

\bibitem{Vyas}
R.G. Vyas, A note on functions of $p(n)$-$\Lambda$-Bounded Variation, J. Indian Math. Soc. (N.S.) 78 (1-4) (2011) 199--204.

\bibitem{Wang3}
H. Wang, Embedding of classes of functions with $\Lambda_{\varphi}$-bounded variation into generalized Lipschitz classes, Math. Inequal.
Appl. 18 (4) (2015) 1463--1471.

\bibitem{Wt1}
D. Waterman, On  convergence of Fourier series of functions of bounded generalized variation, Studia Math. 44 (1972) 107--117.

\bibitem{Wt2}
D. Waterman, On the summability of Fourier series of functions of $\Lambda$-bounded variation, Studia Math. 55 (1976) 87–-95.

\bibitem{Wt3}
D. Waterman, On $\Lambda$-bounded variation, Studia Math. 57 (1976) 33--45.

\bibitem{Wi}
N. Wiener, The quadratic variation of a function and its Fourier coefficients, Massachusetts J. Math. 3 (1924) 72--94.

\bibitem{Wu}
X. Wu, Embedding of classes of functions with bounded $\Phi$-variation into generalized Lipschitz spaces, Acta Math. Hungar. 150 (1) (2016) 247--257.

\bibitem{You}
L.C. Young, Sur une g\'{e}n\'{e}ralisation de la notion de variation de puissance p-i\`{e}me born\'{e}e au sense de M. Wiener, et sur la convergence des
s\'{e}ries de Fourier, C. R. Math. Acad. Sci. Paris 204 (1937) 470--472.

\end{thebibliography}
\end{document}